\documentclass{amsart}
\usepackage{amssymb,amsmath,amsfonts,epsfig,latexsym,amsthm}

\newtheorem{thm}{Theorem}[section]
\newtheorem{defn}[thm]{Definition}
\newtheorem{prop}[thm]{Proposition}
\newtheorem{lem}[thm]{Lemma}
\newtheorem{cor}[thm]{Corollary}

\newtheorem*{KCT}{Klyachko's Classification Theorem}
\newtheorem*{KCC}{Klyachko's Compatibility Condition}

\theoremstyle{definition}
\newtheorem{remk}[thm]{Remark}

\def\A{\ensuremath{\mathbf{A}}}
\def\C{\ensuremath{\mathbf{C}}}
\def\P{\ensuremath{\mathbf{P}}}
\def\Q{\ensuremath{\mathbf{Q}}} 
\def\R{\ensuremath{\mathbf{R}}}
\def\Z{\ensuremath{\mathbf{Z}}}

\def\E{\mathcal{E}}
\def\F{\mathcal{F}}
\def\L{\mathcal{L}}
\def\M{\mathcal{M}}

\def\X{\mathfrak{X}}

\def\MM{\mathbb{M}}

\def\O{\ensuremath{\mathcal{O}}}

\def\Sch{\mathbf{Sch}}
\def\Sets{\mathbf{Sets}}

\def\bu{\mathbf{u}}

\def\FFl{\mathcal{F}\ell}
\def\Fl{F\ell}

\def\PP{P\!P}

\def\fr{{f\!r}}

\DeclareMathOperator{\divisor}{div}
\DeclareMathOperator{\Hom}{Hom}

\DeclareMathOperator{\Spec}{Spec}

\def\<{\ensuremath{\langle}}
\def\>{\ensuremath{\rangle}}

\begin{document}

\title{Moduli of toric vector bundles}

\subjclass[2000]{}
    
\author{Sam Payne}
\address{Stanford University, Mathematics, Bldg. 380, Stanford, CA 94305}
\email{spayne@stanford.edu}
\thanks{Supported by a Graduate Research Fellowship from the NSF and by the Clay Mathematics Institute.  Part of this work was done during a visit to the Institut Mittag-Leffler (Djursholm, Sweden).}

\begin{abstract}
We give a presentation of the moduli stack of toric vector bundles with fixed equivariant total Chern class as a quotient of a fine moduli scheme of framed bundles by a linear group action.  This fine moduli scheme is described explicitly as a locally closed subscheme of a product of partial flag varieties cut out by combinatorially specified rank conditions.  We use this description to show that the moduli of rank three toric vector bundles satisfy Murphy's Law, in the sense of Vakil.  The preliminary sections of the paper give a self-contained introduction to Klyachko's classification of toric vector bundles.  
\end{abstract}

\maketitle

\tableofcontents

\section{Introduction}

In general, moduli of vector bundles are not locally of finite type, even for vector bundles with fixed Chern class on projective space.  One may get around this complication by imposing a stability condition and constructing a moduli space of stable or semistable bundles.  However, if $T \subset PGL_{n+1}$ is a maximal torus, then there does exist a locally finite type moduli stack of $T$-equivariant vector bundles on $\P^n$.  Similarly, there is a locally finite type moduli stack of toric vector bundles on an arbitrary toric variety.  Furthermore, each moduli stack of toric vector bundles with fixed rank and equivariant total Chern class on a toric variety is a global quotient of a quasi-projective scheme by a linear group action, and is canonically defined over $\Spec \Z$.  These properties may make toric vector bundles a convenient testing ground for general questions about vector bundles and their moduli.  The main purpose of this paper is to present a construction of these moduli spaces and to study some of their basic geometric properties.

Let $X = X(\Delta)$ be a toric variety over a field $k$.  A toric vector bundle on $X$ is a vector bundle $\E$ together with a $T$-action on $\E$ compatible with the action on $X$.  Klyachko elegantly classified toric vector bundles on $X$ in terms of finite dimensional $k$-vector spaces $E$ with collections of $\Z$-graded filtrations $\{E^\rho(i)\}$, indexed by the rays $\rho \in \Delta$, that satisfy a compatibility condition \cite{Klyachko90}.  This classification was a major advance in the theory but did not immediately lead to a satisfactory understanding of families of toric vector bundles; in the same paper, Klyachko gave examples, for $k = \C$, of algebraically varying families of filtration $\{E^\rho(i)_t \}$ that satisfy the compatibility condition for every $t \in \C$, but the corresponding toric vector bundles do not all have the same Chern classes.  In particular, a flat family of filtrations satisfying Klyachko's compatibility condition at every geometric point does not necessarily come from a flat family of toric vector bundles.

However, Klyachko's classification did give a very clear understanding of individual toric vector bundles and an elegant means of rigidifying them.  If one ``frames" a toric vector bundle by fixing an isomorphism of one fiber with a trivial vector space, then the toric vector bundle has no nontrivial automorphisms that are compatible with the framing.  The remaining problem of understanding which families of filtrations come from flat families of toric vector bundles motivated recent work on the equivariant Chow cohomology of toric varieties \cite{chowcohom}, which clarified what it means to fix the Chern classes in Klyachko's classification, as follows.  

Let $M$ be the character lattice of $T$, which is canonically identified with the group of integral linear functions on the vector space containing $\Delta$.  For any cone $\sigma \in \Delta$, let $M_\sigma = M / (\sigma^\perp \cap M)$ be the group of integral linear functions on $\sigma$.  Klyachko showed that the filtrations $\{E^\rho(i)\}$ coming from a toric vector bundle satisfy the following compatibility condition.  Let $v_\rho$ denote the primitive generator of a ray $\rho \in \Delta$.  

\begin{KCC} For each cone $\sigma \in \Delta$, there is a decomposition $E = \displaystyle{\bigoplus_{[u] \in M_\sigma}} E_{[u]}$ such that
\[
E^\rho(i) = \sum_{[u](v_\rho) \, \geq \, i} E_{[u]},
\]
 for every $\rho \preceq \sigma$ and $i \in \Z$.
\end{KCC}

\noindent  The splitting $E = \bigoplus E_{[u]}$ is not necessarily unique, but the dimension of $E_{[u]}$ is independent of all choices.  In particular, the multiset $\bu(\sigma)$ in which $[u]$ appears with multiplicity equal to $\dim E_{[u]}$ is well-defined.  These multisets are always compatible with the fan structure in the sense that, if we write $\bu(\sigma)|_\tau$ for the multiset of restrictions to $\tau$ of functions in $\bu(\sigma)$, then
\[
\bu(\tau) = \bu(\sigma)|_\tau,
\]
whenever $\tau \preceq \sigma$.  We write $\Psi_\E = \{ \bu(\sigma) \}_{\sigma \in \Delta}$ for the collection of multisets of linear functions determined by a toric vector bundle $\E$ on $X$.  

In Section \ref{chern classes}, we show that two toric vector bundles $\E$ and $\E'$ have the same equivariant Chern classes if and only if $\Psi_\E = \Psi_{\E'}$.  It follows that if $\{E^\rho(i)_t \}$ is the family of collections of filtrations associated to a flat family of toric vector bundles over a connected base, then the multisets of linear functions $\{ \bu(\sigma)_t \}$ are independent of $t$.  Our main technical result is, roughly, the converse.

\begin{thm}
Let $\Psi = \{ \bu(\sigma) \}_{\sigma \in \Delta}$ be a collection of compatible multisets of linear functions.  Then there is a quasiprojective fine moduli scheme $\M_\Psi^\fr$ for framed toric vector bundles $\E$ on $X$ such that $\Psi_\E = \Psi$.
\end{thm}

\noindent The moduli space $\M^\fr_\Psi$ is canonically defined over $\Spec \Z$ and has a natural embedding as a locally closed subscheme of a product of partial flag varieties which takes a vector bundle to the collection of flags of subspaces appearing in the filtrations $\{E^\rho(i)\}$.  Its image is cut out by rank conditions on the intersections of the various subspaces appearing in the partial flags, and these rank conditions may be expressed combinatorially in terms of $\Psi$.  See Section \ref{rank conditions} for details.  The moduli stack of unframed toric vector bundles $\E$ such that $\Psi_\E = \Psi$ is the quotient of $\M_\Psi^\fr$ by a natural $GL$ action (Section \ref{families}).  In the final section, we apply this construction to show that the moduli of rank three toric vector bundles satisfies Murphy's Law in the sense of Vakil \cite{Vakil06}; in other words, these moduli spaces contain every singularity type of finite type over $\Spec \Z$.  

\begin{remk}
Vakil has shown that many moduli spaces satisfy Murphy's Law, including Hilbert schemes of curves in projective space, moduli spaces of smooth surfaces, and moduli spaces of stable sheaves, but this seems to be the first example of a moduli space of vector bundles with such singularities.
\end{remk}

\begin{remk}
A collection of multisets of linear functions $\{ \bu(\sigma) \}_{\sigma \in \Delta}$ such that $\bu(\tau) = \bu(\sigma)|_\tau$ whenever $\tau \preceq \sigma$ may be conveniently encoded as a piecewise linear function on a ``branched cover" of the fan $\Delta$.  This point of view seems to be useful for studying the set of all vector bundles of fixed rank on a specific toric variety, and is developed in \cite{branchedcovers}.
\end{remk}

\noindent \textbf{Acknowledgments.} I am grateful to B.~Conrad, W.~Fulton, J.~Harris, M.~Hering, A.~Knutson, R.~Lazarsfeld, M.~Musta\c{t}\v{a}, M.~Perling, N.~Proudfoot, K.~Smith, R.~Vakil, and the referee for helpful discussions and suggestions.

\section{Preliminaries} \label{prelim}

Here we give a self-contained exposition of the basic properties of toric vector bundles and their classification.  All of the results in this section appear elsewhere in the literature, some in several places, but since there seems to be no single suitable reference for this material, we have included proofs for the convenience of the reader.  Our presentation draws on \cite{Kaneyama75}, \cite{Klyachko91}, \cite{KnutsonSharpe98}, and \cite{Perling04}, as well as Klyachko's seminal paper \cite{Klyachko90}.  We have attempted to simplify the original proofs where possible.  A brief sketch of the historical development of this subject is included at the end of the section.

\subsection{Notation} \label{notation}

We work over a field $k$.

Let $N$ be a lattice, with $M = \Hom(N, \Z)$ the dual lattice and $T = \Spec k[M]$ the torus with character lattice $M$.  Let $\sigma$ be a convex rational polyhedral cone in $N_\R = N \otimes_\Z \R$, and let
\[
M_\sigma = M / (\sigma^\perp \cap M)
\]
be the group of integral linear functions on $\sigma$.  The map $[u] \mapsto \O(\divisor \chi^{u})$, where $u \in M$ is any lift of $[u] \in M_\sigma$, gives an isomorphism from $M_\sigma$ to the group of toric line bundles on the affine toric variety $U_\sigma$.  In particular, the underlying line bundle of any toric line bundle on $U_\sigma$ is trivial.  We write $\L_{[u]}$ for the toric line bundle on $U_\sigma$ corresponding to a linear function $[u]$.

Let $\Delta$ be a fan in $N_\R$.  We refer to the one-dimensional cones of $\Delta$ as ``rays", and write $v_\rho$ for the primitive generator of a ray $\rho \in \Delta$.  Let $\E$ be a toric vector  bundle on $X$.  Then $T(k)$ acts algebraically on $\Gamma(X, \E)$, the $k$-vector space of global sections, with $t \in T$ acting by
\[
(t \cdot s)(x) = t \, (s \, (t^{-1} x)),
\]
for $s \in \Gamma(X, \E)$ and $x \in X$.  This action induces a decomposition into $T$-eigenspaces
\[
\Gamma(X, \E) = \bigoplus_{u \in M} \Gamma(X, \E)_u,
\]
where $t \cdot s = \chi^u(t) s$ for $t \in T$ and $s \in \Gamma(X, \E)_u$.

\begin{remk}
With this choice of notation, a regular $T$-eigenfunction $\chi^u$ is an element of $\Gamma(X, \O_X)_{-u}$.  Other natural choices are possible, and sign conventions vary widely in the literature.
\end{remk}

\noindent Two sections in $\Gamma(X, \E)_u$ that agree at the identity $x_0$ in the dense torus $T$ must agree on $T$, and hence on all of $X$.  Therefore, writing $E = \E_{x_0}$, evaluation at $x_0$ gives an injection $\Gamma(X, \E)_u \hookrightarrow E$, for each $u \in M$.  We write
\[
E^\sigma_u \subset E
\]
for the image of $\Gamma(U_\sigma, \E)_u$ under evaluation at $x_0$.

For $u' \in (\sigma^\vee \cap M)$, multiplication by $\chi^{u'}$ gives a natural map from $\Gamma(U_\sigma, \E)_u$ to  $\Gamma(U_\sigma, \E)_{u-u'}$.  
This multiplication map commutes with evaluation at $x_0$, inducing an inclusion
\[
E^\sigma_u \subset E^\sigma_{u-u'}.
\]
For  $u' \perp \sigma$, these maps are isomorphisms, so $E^\sigma_u$ depends only on the class $[u] \in M_\sigma$.  For a ray $\rho \in \Delta$ we write
\[
E^\rho(i) = E^\rho_u,
\]
for any $u$ such that $\< u, v_\rho \> = i$.  Then we have a decreasing filtration of $E$
\[
\cdots \supset E^\rho(i-1) \supset E^\rho(i) \supset E^\rho(i+1) \supset \cdots.
\]

\vspace{0 pt}

\subsection{Vector bundles on affine toric varieties}

Gubeladze famously proved that an arbitrary vector bundle on an affine toric variety is trivial \cite{Gubeladze87}, generalizing the Quillen-Suslin Theorem that projective modules over polynomial rings are free.  The analogous property of toric vector bundles is much easier to prove.  
\begin{prop} 
Every toric vector bundle on an affine toric variety splits equivariantly as a sum of toric line bundles whose underlying line bundles are trivial.
\end{prop}

\begin{proof}
Let $\E$ be a toric vector bundle on $U_\sigma$, and let $x_\sigma$ be a point in the minimal $T$-orbit $O_\sigma \subset U_\sigma$.  Choose $T$-eigensections $s_1, \ldots, s_r$ such that $\{ s_i (x_\sigma)\}$ is a basis for the fiber $\E_{x_\sigma}$.  The set of points $x \in U_\sigma$ such that $\{ s_i(x)\}$ is not a basis for $\E_x$ is closed, $T$-invariant, and does not contain $O_\sigma$, and hence is empty.  Therefore, $\E$ splits equivariantly as the sum of the line subbundles spanned by the $s_i$.
\end{proof}

\begin{cor} \label{splitting}
There is a natural bijection between finite multisets $\bu(\sigma) \subset M_\sigma$ and isomorphism classes of toric vector bundles on $U_\sigma$, given by
\[
\bu(\sigma) \mapsto \bigoplus_{[u] \in \bu(\sigma)} \L_{[u]}.
\]
\end{cor}

\begin{remk}
If $\E$ is a toric vector bundle on $U_\sigma$, then the multiset $\bu(\sigma) \subset M_\sigma$ corresponding to the isomorphism class of $\E$ is reflected in the action of the stabilizer  of a point in the minimal $T$-orbit, as follows.  The stabilizer of a point $x_\sigma$ in the minimal orbit $O_\sigma$ is the subtorus $T_\sigma \subset T$ whose character lattice is $M_\sigma$.  Then the multiplicity of $[u]$ in $\bu(\sigma)$ is the dimension of the $\chi^{[u]}$ isotypical component of $\E_{x_\sigma}$, considered as a representation of $T_\sigma$.
\end{remk}

\subsection{Klyachko's classification of toric vector bundles}

Let $E$ and $F$ be $k$-vector spaces with collections of decreasing filtrations $\{E^\rho(i)\}$ and $\{F^\rho(i)\}$, respectively, indexed by the rays $\rho$ in $\Delta$.  A morphism of vector spaces with filtrations is a linear map
\[
\varphi: E \rightarrow F
\]
that takes $E^\rho(i)$ into $F^\rho(i)$ for every $\rho \in \Delta$ and $i \in \Z$.  

\begin{KCT}
The category of toric vector bundles on $X(\Delta)$ is naturally equivalent to the category of finite-dimensional $k$-vector spaces $E$ with collections of decreasing filtrations $\{E^\rho(i)\}$ indexed by the rays of $\Delta$, satisfying the following compatibility condition.

For each cone $\sigma \in \Delta$, there is a decomposition $E = \displaystyle{\bigoplus_{[u] \in M_\sigma}} E_{[u]}$ such that
\[
E^\rho(i) = \sum_{[u](v_\rho) \, \geq \, i} E_{[u]},
\]
 for every $\rho \preceq \sigma$ and $i \in \Z$.
\end{KCT}

\begin{proof}
The equivalence is given, in one direction, by associating to a toric vector bundle $\E$ its fiber $E = \E_{x_0}$ over the identity in the dense torus, together with the filtrations $E^\rho(i)$ defined in Section~\ref{notation}.  A morphism of toric vector bundles $f: \E \rightarrow \F$ maps fibers linearly to fibers, and hence induces a linear map $\varphi_f: E \rightarrow F$.  Furthermore, $f$ respects the eigenspace decompositions of the modules of sections, and in particular takes $\Gamma(U_\rho, \E)_u$ into $\Gamma(U_\rho, \F)_u$ for each ray $\rho$ in $\Delta$ and $u \in M$.  Hence $\varphi_f$ takes $E^\rho(i)$  into $F^\rho(i)$ for every $i \in \Z$, as required.

In the other direction, given a vector space with filtrations $(E, \{ E^\rho(i) \})$ satisfying the compatibility condition, let
\[
E^\sigma_u = \bigcap_{\rho \preceq \sigma} E^\rho(\<u, v_\rho\>),
\]
for $\sigma \in \Delta$ and $u \in M$, and define
\[
E^\sigma = \bigoplus_{u \in M} E^\sigma_u.
\]
Then $E^\sigma$ has a natural $k[U_\sigma]$-module structure, where multiplication by $\chi^{u'}$, for $u' \in (\sigma^\vee \cap M)$, is the sum of the inclusions $E^\sigma_u \subset E^\sigma_{u - u'}$.  Let $T$ act on $E^\sigma$ such that $E^\sigma_u$ is the $\chi^u$-isotypical part.  Then it is straightforward to check using the compatibility condition that the induced quasicoherent sheaf $\widetilde{E}^\sigma$ on $U_\sigma$ is torically isomorphic to 
\[
\bigoplus_{[u] \in M_\sigma} \L_{[u]} \otimes_k E_{[u]}.
\] 
In particular, $\widetilde{E}^\sigma$ is locally free.  Furthermore, the direct sum decomposition implies that $E^\rho(i) = E$ for $i \ll 0$, from which it follows that the natural inclusions
\[
E^\sigma_{u} \subset E^\tau_{u-u'},
\]
for $\tau \preceq \sigma$ and $u' \perp \tau$, induce toric isomorphisms
\[
\widetilde{E}^\sigma|_{U_\tau} \xrightarrow{\sim} \widetilde{E}^\tau.
\]
It is straightforward to check that these isomorphisms are gluing data, that morphisms of vector spaces with compatible filtrations induce morphisms of toric vector bundles, and that the functor so defined is inverse to the functor $\E \mapsto (E, \{E^\rho(i)\})$, up to natural isomorphisms, giving an equivalence of categories.
\end{proof}

\subsection{Historical background}

Kaneyama gave the first classification of toric vector bundles in 1975 \cite{Kaneyama75}, only a few years after mathematicians began studying toric varieties systematically in \cite{Demazure70}, \cite{KKMS}, and \cite{MiyakeOda75}.  Although stated only for smooth complete toric varieties, Kaneyama's classification can be extended in a straightforward way to arbitrary toric varieties.  The classification involves both combinatorial and linear algebraic data, as well as a nontrivial equivalence relation.  The combinatorial part, modulo the equivalence relation, encodes the splitting type of the vector bundle on each maximal $T$-invariant affine open subvariety.  The linear algebraic part consists of an element of $GL_r$ for each pair of maximal cones, satisfying a cocycle condition and a compatibility condition with the combinatorial data, and encodes the transition functions.  Bertin and Elencwajg used Kaneyama's results and other equivariant methods to study vector bundles of small rank on projective spaces \cite{BertinElencwajg82}, and Kaneyama eventually used this approach to give a complete description of toric vector bundles of rank at most $n$ on $\P^n$ \cite{Kaneyama88}.  Similar results were obtained independently by Behrmann, in his thesis \cite{Behrmann86}.

Klyachko introduced his classification of toric vector bundles in terms of vector spaces with compatible filtrations in the late 1980s; see \cite{Klyachko87}, \cite{Klyachko89}, \cite{Klyachko90}, and \cite{Klyachko90b}.  His original proof of the classification made extensive use of limits of sections over one-parameter subgroups of the dense torus.  Klyachko used this classification most notably to study stability conditions on toric vector bundles, especially on $\P^2$, from which he gained deep insight into the Horn Problem on eigenvalues of sums of Hermitian matrices \cite{Klyachko98}.  Toric vector bundles on $\P^2$ also appear in Penacchio's work on Hodge structures \cite{Penacchio03}.

In a 1991 preprint \cite{Klyachko91}, Klyachko stated an extension of his classification to toric reflexive and torsion free coherent sheaves, with filtrations of the generic fiber constructed from the $T$-eigenspace decompositions of the modules of sections over $T$-invariant affine open subvarieties.  These ideas were pushed a little further, and presented with many helpful expository details, as well as some connections to string theory, by Knutson and Sharpe in \cite{KnutsonSharpe98}.  Further connections between toric vector bundles, string theory, and mirror symmetry are explored in \cite{KnutsonSharpe99} and \cite{LiuYau00}.  

Recently, Perling has systematically developed the approach to toric sheaves through $T$-eigenspace decompositions of modules of sections, giving a classification of arbitrary toric quasicoherent sheaves \cite{Perling04}.  Perling's classification data, which he calls $\Delta$-families of vector spaces, explicitly encode the $T$-eigenspace decompositions of the modules of sections, together with the multiplication maps for regular $T$-eigenfunctions.  He has applied this approach to study moduli of rank two toric vector bundles on toric surfaces, using equivariant resolutions analogous to Euler sequences \cite{Perling04b}, as well as resolutions of more general toric sheaves \cite{Perling05}.
 
Another classification of equivariant quasicoherent sheaves on toric varieties, developed earlier by Cox  \cite{Cox95} and Musta\c{t}\v{a} \cite{Mustata02}, involves fine graded modules over the homogeneous coordinate ring.  Although fine graded modules and $\Delta$-families both use multigraded commutative algebra, they are very different in flavor.  Perling's $\Delta$-families are specifically tailored to toric sheaves and lead to elegant characterizations of toric reflexive and locally free sheaves within this broader framework.   On the other hand, graded modules over homogeneous coordinate rings work equally well for studying arbitrary (not necessarily toric) sheaves, but seem to be less natural for studying reflexive sheaves or vector bundles specifically.  Both approaches are likely to be useful for future research.

\section{Constructing the moduli spaces} \label{moduli}

\subsection{Equivariant Chern classes} \label{chern classes}

Let $\Delta$ be a fan, and let $X = X(\Delta)$.  Recall that the equivariant Chow cohomology ring $A^*_T(X)$ is naturally isomorphic to the ring of integral piecewise polynomial functions
\[
\PP^*(\Delta) = \{ f: |\Delta| \rightarrow \R \, | \, f|_\sigma  \in M_\sigma \mbox{ for all } \sigma \in \Delta \}.
\]
The isomorphism takes a linear function $u \in M$ to the first Chern class of the toric line bundle $\O_X(\divisor \chi^u)$ \cite{chowcohom}.  

Let $\Psi = \{ \bu(\sigma) \}_{\sigma \in \Delta}$ be a collection of multisets of linear functions of size $r$ such that
\[
\bu(\tau) = \bu(\sigma)|_\tau,
\]
for $\tau \preceq \sigma$.  Let $c_i(\Psi)$ be the piecewise polynomial function whose restriction to $\sigma$ is $e_i(\bu(\sigma))$, the $i$-th elementary symmetric function in the multiset of linear functions $\bu(\sigma)$.  We define
\[
c(\Psi) = 1 + c_1(\Psi) + \cdots + c_r(\Psi).
\]

\noindent Recall that we write $\Psi_\E$ for the collection of compatible multisets of linear functions determined by a toric vector bundle $\E$.

\begin{prop}
The equivariant total Chern class of a toric vector bundle $E$ is $c_T(\E) = c(\Psi_\E)$.
\end{prop}

\begin{proof}
Since the restriction of $\E$ to $U_\sigma$ is 
\[
\E|_{U_\sigma} \cong \bigoplus_{[u] \in \bu(\sigma)} \L_{[u]},
\]
and since the first Chern class of $\L_{[u]}$ corresponds to the linear function $[u]$, the proposition follows from the naturality of the isomorphism $A^*_T(X) \cong \PP^*(\Delta)$.
\end{proof}

 \subsection{Flags and rank conditions} \label{rank conditions}
 
 In this section, we interpret the compatibility condition on the filtrations in Klyachko's Classification Theorem in terms of rank conditions on the flags of subspaces appearing in these filtrations, for toric vector bundles with fixed equivariant Chern classes.  This discussion motivates and gives the intuition for the more technical arguments in Section~\ref{families}.
 
Suppose $\E$ is an equivariant vector bundle on $X(\Delta)$ with equivariant total Chern class
\[
c_T(\E) = c(\Psi).
\]
For each ray $\rho$ in $\Delta$, let $J(\rho) \subset \{1, \ldots, r\}$ be the set of dimensions of the nonzero vector spaces appearing in the filtration $E^\rho(i)$ of the fiber $E = \E_{x_0}$.  Let
\[
\Fl(\rho) \in \FFl_{J(\rho)} (E)
\]
be the partial flag in $E$ consisting of exactly these subspaces.  It follows from Klyachko's Classification Theorem that $\E$ is determined up to isomorphism by the data $(\Psi, \{ \Fl(\rho) \})$.  To see this, let $\Lambda_\rho : J(\rho) \rightarrow \Z$ be the function given by
\[
\Lambda_\rho(j) = \max \{ i \ | \ \dim E^\rho(i) = j \},
\]
and let $\Fl(\rho)_j$ be the $j$-dimensional subspace of $E$ appearing in $\Fl(\rho)$, for $j \in J(\rho)$.  The filtrations determining $E$ are then given by
\[
E^\rho(i) = \sum_{\Lambda_\rho(j) \, \geq \, i} \Fl(\rho)_j.
\]

The fixed collection of multisets of linear functions $\Psi$ restricts how the partial flags $\{ \Fl(\rho) \}$ can meet each other, as follows.  Klyachko's compatibility condition says that there is a splitting $E = \displaystyle{\bigoplus_{[u] \in \bu(\sigma)} E_{[u]}}$ such that
\[
E^\rho(i) = \bigoplus_{[u](v_\rho) \, \geq \, i} E_{[u]},
\]
for every ray $\rho \preceq \sigma$.  Therefore, if $\rho_1, \ldots, \rho_s$ are rays of $\sigma$, with $j_\ell \in J(\rho_\ell)$ for $1 \leq \ell \leq s$, then
\[
\dim \bigcap_{\ell=1}^s \Fl(\rho_\ell)_{j_\ell}  \ = \ \# \{ [u] \in \bu(\sigma) \ | \ [u](v_{\rho_\ell}) \geq \Lambda_{\rho_l}(j_\ell) \mbox{ for } 1 \leq \ell \leq s \}.
\]

\vspace{-27 pt} \hfill (*) \vspace{17 pt}

\noindent We refer to (*) as a rank condition because it can be interpreted as specifying the rank of the natural map
\[
E \rightarrow \prod_{\ell=1}^s E/ \Fl(\rho_\ell)_{j_\ell},
\]
since the kernel of this map is exactly $\bigcap_{\ell = 1}^s \Fl(\rho_\ell)_{j_\ell}$.

After a choice of isomorphism $E \cong k^r$, each rank condition (*) can be written as the vanishing of certain polynomials in the Pl\"ucker coordinates of the partial flags $\Fl(\rho)$ plus the nonvanishing of certain other such polynomials.  These rank conditions, for all cones $\sigma \in \Delta$, and all collections of rays $\{ \rho_1, \ldots, \rho_s \}$ of $\sigma$, and all $\{j_\ell\}_{1 \leq \ell \leq s}$, cut out a locally closed subscheme of the product of partial flag varieties
\[
\FFl_{\Psi} = \prod_\rho \FFl_{J(\rho)}(k^r).
\]
We write $\M^\fr_{\Psi}$ for this locally closed subscheme.  In Section \ref{families}, we will show that $\M^\fr_{\Psi}$ is a fine moduli space for framed toric vector bundles on $X$ with equivariant total chern class $c(\Psi)$.  See Definition \ref{framed bundles} and Theorem \ref{framed moduli}.  

\begin{remk}
The schemes $\M_{\Psi}^\fr$ are similar in flavor to the permutation array schemes studied by Billey and Vakil in \cite{BilleyVakil05}, which are also locally closed subschemes of products of flag varieties cut out by combinatorially encoded rank conditions.  One minor difference is that permutation array schemes are cut out by one condition for every possible intersection among the different flags, while the conditions cutting out $\M_{\Psi}^\fr$ are given by only a subset of the possible intersections.
\end{remk}

\begin{defn} \label{framed bundles}
A framed toric vector bundle on $X$ is a toric vector bundle $\E$ together with an isomorphism $\varphi: E \xrightarrow{\sim} k^r$.
\end{defn}

\noindent A morphism of framed toric vector bundles is a morphism of toric vector bundles that is compatible with the framings.

\begin{prop} \label{weak moduli}
The map $\E \mapsto \{ \Fl(\rho) \}$ gives a bijection between the set of isomorphism classes of rank $r$ framed toric vector bundles on $X$ with equivariant total Chern class $c(\Psi)$ and the set of $k$-points of $\M^\fr_{\Psi}$.
\end{prop}

\begin{proof}
The injectivity of this map follows from Klyachko's Classification Theorem.  Surjectivity is an immediate consequence of the following lemma.
\end{proof}

\begin{lem} \label{rank lemma}
Let $\bu(\sigma) \subset M_\sigma$ be a multiset, and let $(E, \{E^\rho(i)\})$ be a vector space with a collection of filtrations indexed by the rays of $\sigma$ such that, for any collection of integers $\{ i_\rho \}$,
\[
\dim \Bigl( \, \bigcap_\rho E^\rho(i_\rho) \Bigr) \, = \, \# \bigl\{ [u] \in \bu(\sigma) \ | \ [u](v_\rho) \geq i_\rho \mbox{ for all } \rho \bigr\}.
\]
Then there is a splitting $E = \bigoplus_{[u] \in \bu(\sigma)} E_{[u]}$ such that
\[
E^\rho(i) = \sum_{[u](v_\rho) \, \geq \, i} E_{[u]},
\]
for all $\rho$ and all $i$.
\end{lem}

\begin{proof}
Consider the partial order  on $\bu(\sigma)$ where $[u] \leq [u']$ if and only if $[u]-[u']$ is nonnegative on $\sigma$.  We choose $E_{[u]}$ inductively by moving up through $\bu(\sigma)$ with respect to this partial order.  If $[u]$ is minimal, then $\dim \bigcap_\rho E^\rho([u](v_\rho))$ is the multiplicity of $[u]$ in $\bu(\sigma)$, and we set
\[
E_{[u]} = \bigcap_\rho E^\rho \bigl( [u](v_\rho) \bigr).
\]
If $[u]$ is not minimal, then
\[
\dim \Bigl( \, \bigcap_\rho E^\rho \bigl( [u](v_\rho) \bigr) \Bigr) \, = \, \# \bigl\{ [u'] \in \bu(\sigma) \ | \ [u'] \leq [u] \bigr\}.
\]
Assuming, by induction, that $E_{[u']}$ is fixed for each $[u'] \leq [u]$, we choose $E_{[u]}$ to be a subspace of $\bigcap_\rho E^\rho([u](v_\rho))$ complementary to $\sum_{[u'] < [u]} E_{[u']}$.  For any such choice of $\{ E_{[u]} \}$, $E$ decomposes as a direct sum $E = \bigoplus E_{[u]}$, and it is straightforward to check that
\[
E^\rho(i) = \sum_{[u](v_\rho) \, \geq \, i} E_{[u]},
\]
for any ray $\rho$ and any integer $i$.
\end{proof}

\begin{cor} \label{weak moduli quotient}
There is a natural bijection between the set of isomorphism classes of rank $r$ toric vector bundles on $X$ with equivariant total Chern class $c(\Psi)$ and the set of $GL_r(k)$-orbits of $\M^\fr_{\Psi}(k)$.
\end{cor}

\subsection{Families of toric vector bundles} \label{families}

In this section we study families of framed and unframed toric vector bundles parametrized by a scheme $S$ and the related moduli problems.

As in Section~\ref{rank conditions}, we fix a fan $\Delta$ and a collection $\Psi = \{ \bu(\sigma) \}_{\sigma \in \Delta}$ of multisets of linear functions of size $r$  such that $\bu(\tau) = \bu(\sigma)|_\tau$, for $\tau \preceq \sigma$.  By Proposition~\ref{weak moduli}, isomorphism classes of framed rank $r$ toric vector bundles on $X = X(\Delta)$ with equivariant total Chern class $c(\Psi)$ are naturally parametrized by the $k$-points of the locally closed subscheme
\[
\M^\fr_{\Psi} \subset \FFl_{\Psi}
\]
cut out by the rank conditions (*).

The main result of this section is Theorem~\ref{framed moduli}, which says that the isomorphism classes of $S$-families of framed rank $r$ toric vector bundles on $X$ with equivariant total Chern class $c(\Psi)$ are naturally parametrized by the $S$-points of $\M^\fr_{\Psi}$, for any scheme $S$.  The key step is to prove a relative version of Klyachko's Classification Theorem, with the splittings in the compatibility condition replaced by rank conditions, using a relative version of the isotypical decomposition of the modules of sections over $T$-invariant open subvarieties \cite[I.4.7.3]{SGA3.1}.  See Proposition~\ref{relative KCT}.

From this point on, we work over $\Spec \Z$, instead of over a field $k$.  This is natural because toric varieties, flag varieties, and the locally closed subschemes $\M^\fr_{\Psi}$ are all canonically defined over $\Spec \Z$ and, once we are working in the relative setting, this additional level of generality does not create any new difficulties or complications.  We write $\Sch$ for the category of schemes of finite type over $\Spec \Z$.

\vspace{5 pt}

Let $S$ be a scheme of finite type over $\Spec \Z$, and let $T_S$ be the relative torus $T \times S$.

\begin{defn}
An $S$-family of  toric vector bundles on $X$ is a vector bundle $\E$ on $X \times S$ with an action of the relative torus $T_S$ compatible with the action on $X \times S$.
\end{defn}

\noindent We say that $\E$ is a family of toric vector bundles with total Chern class $c(\Psi)$ if 
\[
c_T(\E|_{X \times s}) = c(\Psi),
\]
for every geometric point $s$ of $S$.  

A morphism of $S$-families of toric vector bundles on $X$ is a morphism of vector bundles over $X \times S$ that is compatible with the actions of $T_S$.  We write $\MM_{\Psi} : \Sch \rightarrow \Sets$ for the moduli functor
\[
\MM_{\Psi}(S) = \left\{ \begin{array}{c} \mbox{isomorphism classes of $S$-families of rank $r$ toric vector} \\\mbox{bundles on $X$ with equivariant total Chern class $c(\Psi)$} \end{array} \right\}.
\]
One cannot hope for the moduli functor $\MM_{\Psi}$ to be representable by a scheme in general, but in many cases it has a coarse moduli scheme.  This coarse moduli scheme, when it exists, may be constructed as a quotient of a fine moduli scheme for framed toric vector bundles.

\begin{defn}
 An $S$-family of framed toric vector bundles on $X$ is a family of toric vector bundles with an isomorphism $\varphi: \E|_{x_0 \times S} \rightarrow \O_S^{\oplus r}$.
 \end{defn}

\noindent A morphism of $S$-families of framed toric vector bundles is a morphism of $S$-families of toric vector bundles compatible with the framings.  We write $\MM^\fr_{\Psi}: \Sch \rightarrow \Sets$ for the moduli functor
\[
\MM^\fr_{\Psi} (S) = \left\{ \begin{array}{c} \mbox{isomorphism classes of $S$ families of framed rank $r$ toric} \\ \mbox{vector bundles with equivariant total Chern class $c(\Psi)$} \end{array} \right\}.
\]

\begin{thm} \label{framed moduli}
There is a natural isomorphism of functors
\[
\MM^\fr_{\Psi} \ \cong \ \Hom \bigl( \, \underline{\ \ \ } \, , \, \M^\fr_{\Psi} \bigr).
\]
\end{thm}

\noindent  Before giving the proof, we use Theorem \ref{framed moduli} to describe the coarse moduli scheme for $\MM_{\Psi}$, when it exists.

Each of the rank conditions (*) cutting out $\M^\fr_{\Psi}$ is invariant under the diagonal action of $GL_r$ on the product of partial flag varieties $\FFl_{\Psi}$.  Therefore, the diagonal action of $GL_r$ restricts to an action on $\M^\fr_{\Psi}$. 

\begin{prop}\label{coarse moduli}
A scheme $\M$ is a coarse moduli space for $\MM_{\Psi}$ if and only if it is a categorical quotient of $\M^\fr_{\Psi}$ for the action of $GL_r$, and, for every algebraically closed field $k$, the set of $GL_r(k)$-orbits of $\M^\fr_{\Psi}(k)$ maps bijectively onto $\M(k)$. 
\end{prop}

\begin{proof}
This follows from standard arguments used in the construction of moduli of coherent sheaves on projective varieties.  See, for instance, \cite[Chapter 4]{HuybrechtsLehn97}.
\end{proof}

The center $\mathbb{G}_m$ of $GL_r$, consisting of invertible scalar multiples of the identity, acts trivially on $\FFl_{\Psi}$, so the quotient $PGL_r$ acts naturally on $\FFl_{\Psi}$ and on $\M^\fr_{\Psi}$.  There are many cases where $PGL_r$ acts freely on $\M^\fr_{\Psi}$ (see Section \ref{Murphy}).

\begin{cor}
If $PGL_r$ acts freely on $\M^\fr_{\Psi}$, then there is a geometric quotient $\M^\fr_{\Psi} /\!/ GL_r$, and it is a coarse moduli scheme for $\MM_{\Psi}$.
\end{cor}

\begin{proof}
If $PGL_r$ acts freely on $\M^\fr_{\Psi}$, then the action of $GL_r$ is closed.  The existence of the geometric quotient $\M^\fr_{\Psi} /\!/ GL_r$ then follows from Geometric Invariant Theory for reductive groups \cite[Chapter 1]{GIT}, and it is a coarse moduli scheme for $\MM_{\Psi}$ by Proposition \ref{coarse moduli}.
\end{proof}

\begin{remk}
For those who prefer such language, the stack $\X_{\Psi}$, whose objects over $S \in \Sch$ are $S$-families of toric vector bundles, with morphisms given by Cartesian squares, is isomorphic to the quotient stack $[\M^\fr_{\Psi} / GL_r]$.  The map from $\X_{\Psi}$ to $[\M^\fr_{\Psi} / GL_r]$ takes an $S$-family of toric vector bundles to the pullback of the universal frame bundle under the induced map $S \rightarrow \M^\fr_{\Psi}$.  The map in the other direction takes a $GL_r$-torsor $Y$ over $S$ with an equivariant map $Y \rightarrow \M^\fr_{\Psi}$ to the quotient of the induced $Y$-family of framed toric vector bundles by the natural $GL_r$-action.  For a short and gentle introduction to stacks that includes the basic technical details arising in this example, see \cite{Fantechi01}.
\end{remk}

It remains to prove Theorem \ref{framed moduli}.  The key is the following relative version of Klyachko's Classification Theorem.  Note that, in the relative setting, one cannot ask for the splittings that appear in Klyachko's compatibility condition.  Instead, we use rank conditions, as in Lemma \ref{rank lemma}.

\begin{prop}\label{relative KCT}
The category of $S$-families of toric vector bundles on $X$ is naturally equivalent to the category of vector bundles $E$ on $S$ with collections of decreasing filtrations $\{ E^\rho(i) \}$ indexed by the rays of $\Delta$, satisfying the following rank condition.

For each cone $\sigma \in \Delta$, there is a multiset $\bu(\sigma) \subset M_\sigma$ such that, for any rays $\rho_1, \ldots, \rho_s$ of $\sigma$ and integers $i_1, \ldots, i_s$, $E^{\rho_1}(i_1) \cap \cdots \cap E^{\rho_s}(i_s)$ is a vector bundle of rank equal to 
\[
\# \{ [u] \in \bu(\sigma) \ | \ [u](v_{\rho_j}) \geq i_j \mbox{ for } 1 \leq j \leq s \}.
\]
\end{prop}

\begin{proof}
The equivalence is given, in one direction, by associating to an $S$-family of toric vector bundles $\E$, the vector bundle $E = \E|_{x_0 \times S}$, with filtrations $E^\rho(i)$ defined as follows.  By \cite[I.4.7.3]{SGA3.1}, for each $\sigma \in \Delta$, $\E|_{U_\sigma \times S}$ splits as a direct sum of $\O_S$-modules
\[
\E|_{U_\sigma \times S} \cong \bigoplus_{u \in M} E^\sigma_u,
\]
where $E^\sigma_u$ is $\chi^u$-isotypical for the action of $T_S$.  As in the proof of Klyachko's Classification Theorem, restriction to $x_0 \times S$ gives injections $E^\sigma_u \subset E$ whose image depends only on the class $[u]$ in $M_\sigma$.  Define $E^\rho(i) \subset E$ to be the image of $E^\rho_u$ for any $u$ such that $\< u, v_\rho \> = i$.  Each $E^\rho(i)$ is $\O_S$-coherent and $\O_S$-flat, and hence is a vector bundle.  The rank condition follows from Klyachko's Classification Theorem and Lemma \ref{rank lemma} applied to $\E|_{X \times s}$, for any geometric point $s \in S$.  A morphism of $S$-families of toric vector bundles respects the $T_S$ actions, so restricting to $x_0 \times S$ gives a morphism of vector bundles with filtrations.

In the other direction, given a vector bundle with filtrations $(E, \{ E^\rho(i) \})$ satisfying the rank condition, let $E^\sigma_u$ be the vector bundle on $S$ given by
\[
E^\sigma_u  = \bigcap_{\rho \preceq \sigma} E^\rho(\<u, v_\rho \>),
\]
for $\sigma \in \Sigma$ and $u \in M$.  Define $E^\sigma = \bigoplus_{u \in M} E^\sigma_u$.  Just as in the proof of Klyachko's Classification Theorem, there is an induced $T_S$-equivariant sheaf $\widetilde{E}^\sigma$ on $U_\sigma \times S$, and the sheaves $\{ \widetilde{E}^\sigma \}_{\sigma \in \Delta}$ glue together to give an $S$-family of toric vector bundles $\E$.  It is straightforward to check that morphisms of vector bundles with filtrations on $S$ induce morphisms of $S$-families of toric vector bundles, and that the functor so defined is inverse to the functor $\E \mapsto (E, \{E^\rho(i) \})$, up to natural isomorphisms, giving an equivalence of categories.
\end{proof}

\begin{proof}[Proof of Theorem \ref{framed moduli}]
Let $\E$ be an $S$-family of rank $r$ toric vector bundles on $X$ with equivariant total Chern class $c(\Psi)$, and let $\varphi: \E|_{x_0 \times S} \xrightarrow{\sim} \O_S^{\oplus r}$ be a framing of $\E$.  For each ray $\rho$ in $\Sigma$, let $J(\rho) \subset \{1, \ldots, r\}$ be the set of ranks of the bundles $E^\rho(i)$ appearing in the filtration associated to $\E$, and let
\[
\Fl(\rho) \in \FFl_{J(\rho)}(\O_S^{\oplus r})
\]
be the flag consisting of exactly these subbundles.  Then it follows from Proposition~\ref{relative KCT} that $\E \mapsto \{ \Fl(\rho) \}$ gives a bijection between $\MM^\fr_{\Psi}(S)$ and the set of collections of partial flags in $\prod_\rho \FFl_{J(\rho)}(\O_S^{\oplus r})$ satisfying the rank conditions, which is canonically identified with $\Hom(S, \M^\fr_{\Psi})$.
\end{proof}

\vspace{0 pt}

\section{Murphy's Law for toric vector bundles} \label{Murphy}

In this section, we study the singularities of the moduli of toric vector bundles.  For rank two bundles, each fine moduli scheme $\M^\Psi_\fr$  is isomorphic to a product of copies of $\P^1$, minus some collection of diagonals, so in particular it is smooth.  If $\M_{\Psi}^\fr$ is contained in the complement of the diagonals $D_{12} \cup D_{23}$ for some numbering of the $\P^1$ factors, then $PGL_2$ acts simply transitively on the first three factors, so there is a geometric quotient $\M_{\Psi}^\fr /\!/ GL_2$ that is the coarse moduli scheme for toric vector bundles on $X$ with equivariant total Chern class $c(\Psi)$, and it is also smooth.

The main result of this section, roughly stated, is that the moduli of rank three toric vector bundles have arbitrarily bad singularities.  Following \cite{Vakil06}, we define a singularity type to be an equivalence class of pointed schemes, where the equivalence relation is generated by setting a pointed scheme $(X, p)$ to be equivalent to another pointed scheme $(Y,q)$ if there is a smooth morphism $f:X \rightarrow Y$ such that $f(p) = q$, and we say that a scheme $\M$ satisfies Murphy's Law if every singularity type of finite type over $\Spec \Z$ is equivalent to $(\M, p)$ for some point $p$ in $\M$.

By the fine moduli scheme of framed rank $r$ toric vector bundles on a given class of toric varieties, we mean the disjoint union over all toric varieties $X(\Delta)$ in this class and all possible equivariant total Chern classes of the fine moduli scheme of rank $r$ framed toric vector bundles on $X(\Delta)$ with that Chern class.  The coarse moduli scheme of rank $r$ toric vector bundles on a class of toric varieties is defined analogously, but with the disjoint union taken over those Chern classes such that the coarse moduli scheme exists.

\begin{thm} \label{Murphy thm}
The following spaces satisfy Murphy's Law.
\begin{enumerate}
\item The fine moduli scheme of framed rank three toric vector bundles on smooth quasi-affine toric varieties.
\item The fine moduli scheme of framed rank three toric vector bundles on $\Q$-factorial quasiprojective  toric fourfolds.
\item The coarse moduli scheme of rank three toric vector bundles on smooth quasi-affine toric varieties.
\item The coarse moduli scheme of rank three toric vector bundles on $\Q$-factorial quasiprojective toric fourfolds.
\end{enumerate}
\end{thm}

Theorem \ref{Murphy thm} is an immediate consequence of the following stronger result, which implies, roughly speaking, that the moduli of toric vector bundles with fixed equivariant Chern class can be chosen to have arbitrarily many components, with arbitrary singularities along the generic point of each component.  The proof of this result is a straightforward application of Mn\"ev's Universality Theorem.

\begin{thm} \label{universality}
Let $Y$ be an affine scheme of finite type over $\Spec \Z$.  Then there exists a fan $\Delta$ and a collection of compatible multisets of linear functions $\Psi$ such that $PGL_3$ acts freely on $\M_{\Psi}^\fr$ and the geometric quotient $\M_{\Psi}^\fr /\!/ GL_3$ is isomorphic to an open subvariety of $Y \times \A^s$ that projects surjectively onto $Y$, for some positive integer $s$.  Furthermore, $\Delta$ may be chosen such that $X(\Delta)$ is smooth and quasi-affine, or such that $X(\Delta)$ is a $\Q$-factorial quasiprojective fourfold.
\end{thm}

\noindent The rough idea of the proof of Theorem \ref{universality} is as follows.  For any collection $\Psi$ of compatible multisets of linear functions, $\M_{\Psi}^\fr$ is a scheme parametrizing collections of points and lines in $\P^2$ with some specified incidence conditions.  By Mn\"ev's Universality Theorem, it suffices to prove that $\Delta$ and $\Psi$ may be chosen such that these incidence relations are arbitrary, which we will show combinatorially.

Given positive integers $d$ and $d'$ and a subset $I \subset \{1, \ldots, d\} \times \{1, \ldots, d' \}$, let $C^{d,d'}_I$ be the locally closed subscheme
\[
C^{d,d'}_I \subset \ \prod_{i=1}^d \P^2 \times \prod_{j=1}^{d'} \check{\P}{^2}
\]
that parametrizes distinct points $(x_1, \ldots, x_d)$ and distinct lines $(\ell_1, \ldots, \ell_{d'})$ in $\P^2$ such that $x_i$ lies on $\ell_j$ if and only if $(i,j)$ is in $I$.

\begin{proof}[Proof of Theorem \ref{universality}]
By Mn\"ev's Universality Theorem, as presented in \cite[Section 1.8]{Lafforgue03}, there exist positive integers $d$ and $d'$ and a subset $I \subset  \{1, \ldots, d\} \times \{1, \ldots, d' \}$ such that $PGL_3$ acts freely on $C^{d,d'}_I$ and such that the geometric quotient $C^{d,d'}_I /\!/ GL_3$ is isomorphic to an open subscheme of $Y \times \A^s$ that projects surjectively on $Y$.  It remains to show that $\Delta$ and $\Psi$ may be chosen such that $\M_{\Psi}^\fr$ is $PGL_3$-equivariantly isomorphic to $C^{d,d'}_I$.

We first prove the theorem with $X(\Delta)$ smooth and quasi-affine.  Let $N = \Z^{d + d'}$, and let $\Delta$ be the fan in $N_\R$ whose maximal cones are
\[
\sigma_{ij} = \<e_i, e_j\>,
\]
for $1 \leq i,j \leq d + d'$, for $i \neq j$.  Then $X(\Delta)$ is the complement in $\A^{d+d'}$ of the union of its codimension three coordinate subspaces; in particular $X(\Delta)$ is smooth and quasi-affine.  To specify $\Psi$, we give multisets $\bu(\sigma_{ij}) \subset M_{\sigma_{ij}}$ for the maximal cones in $\Delta$ such that, for any $\tau \in \Delta$ the multiset of linear functions $\bu(\sigma_{ij})|_\tau$ is independent of the choice of $\sigma_{ij}$ containing $\tau$.  Define $\bu(\sigma_{ij})$ as follows.
\[
\bu(\sigma_{ij}) = \left\{ \begin{array}{ll}
						\{ 0, e_i^*, e_j^* \} & \mbox{ for } i \leq d \mbox{ and } j \leq d, \\
						\{0, e_j^*, e_i^* + e_j^* \} & \mbox{ for } i \leq d, j > d, \mbox{ and } (i, j-d) \in I, \\
						\{e_i^*, e_j^*, e_j^*\} & \mbox{ for } i \leq d, j > d, \mbox{ and } (i, j-d) \not \in I, \\
						\{e_i^*, e_j^*, e_i^* + e_j^*\} & \mbox{ for } i > d \mbox{ and } j > d.
\end{array} \right.
\]
Then $\FFl_{J(\rho_i)} = \P^2$ for $i \leq d$ and $\FFl_{J(\rho_j)} = \check{\P}{^2}$ for $j > d$, and it is straightforward to check that the rank conditions cutting out the locally closed subscheme
\[
\M_{\Psi}^\fr \subset \prod_{i=1}^d \P^2 \times \prod_{j=d+1}^{d+d'} \check{\P}{^2}
\]
are exactly the conditions that $x_i$ lies on $\ell_j$ if and only if $(i,j-d)$ is in $I$.  This proves the result with $X(\Delta)$ smooth and quasi-affine.

Now consider $N' = \Z^4$.  Let $P \subset N'_\R$ be an embedding of a cyclic four-dimensional polytope with $d + d'$ rational vertices that contains $0$ in its interior. Let $v_1, \ldots, v_{d+d'}$ be the primitive generators of the rays through the vertices of $P$, and let $\Delta'$ be the fan in $N'_\R$ whose maximal cones are $\sigma'_{ij} = \< v_i, v_j\>$.  Then $X(\Delta')$ is $\Q$-factorial, quasiprojective, and four-dimensional.

The projection $\Z^{d+d'} \rightarrow N'$ given by $e_i \mapsto v_i$ maps the fan $\Delta$ onto $\Delta'$, and induces a bijection on the underlying sets $|\Delta| \rightarrow |\Delta'|$.  Then each element of $\bu(\sigma_{ij}))$ is the pullback of some element of $(M_{\sigma_{ij}})_\Q$.  Therefore, for a suitable integer $m$, the collection $m\Psi$ of multisets of linear functions in $\bu(\sigma_ij)$ multiplied by $m$, is the pullback of a collection $\Psi'$ of compatible multisets of linear functions on $\Delta'$.  Then the moduli space $\M_{m \Psi'}^\fr$ of toric vector bundles on $X(\Delta')$ lives in exactly the same product of partial flag varieties as $\M_{\Psi}^\fr$ and is cut out by exactly the same rank conditions, so $\M_{m\Psi'}^\fr$ is also naturally isomorphic to $C^{d,d'}_I$, as required.
\end{proof}

\begin{remk}
It follows easily from the proof of Theorem \ref{universality} that similar universality statements and Murphy's Law statements hold for moduli of rank $r$ toric vector bundles for any $r \geq 3$.
\end{remk}

\begin{remk}
It is not known whether the moduli of toric vector bundles on projective toric varieties satisfy Murphy's Law.
\end{remk}

\begin{remk}
Although this paper deals mainly with separated toric varieties, the main results also hold for toric vector bundles on toric prevarieties when rephrased in terms of multifans.  In particular, Klyachko's Classification Theorem works equally well for toric prevarieties.  The equivariant Chow cohomology ring of a toric prevariety is canonically isomorphic to the ring of piecewise-polynomial functions on the corresponding multifan \cite[Section 4]{chowcohom}, and the moduli of framed toric vector bundles with fixed equivariant Chern class may be constructed just as for separated toric varieties, with essentially identical proofs.  The fans used in the proof of Theorem \ref{universality} are two-dimensional, and a general projection to $\Z^2$ gives a multifan corresponding to a two-dimensional toric prevariety.  In particular, $X(\Delta)$ may be chosen to be a possibly nonseparated toric surface in Theorem \ref{universality}, and Murphy's Law holds for the moduli of rank three toric vector bundles on two-dimensional toric prevarieties.
\end{remk}

\bibliography{math}
\bibliographystyle{amsalpha}

\end{document}